\documentclass{article}
\usepackage{graphicx} 

\usepackage{amsmath, epsfig, cite}
\usepackage{amssymb}
\usepackage{amsfonts}
\usepackage{latexsym}
\usepackage{float}
\usepackage{color}
\usepackage{mathrsfs}

\newtheorem{thm}{Theorem}[section]

\newtheorem{cor}[thm]{Corollary}

\newtheorem{lem}[thm]{Lemma}

\newcommand{\pf}{\noindent{\it Proof} }

\numberwithin{equation}{section}

\newcommand{\qed}{{\hfill$\square$}\medskip}

\begin{document}
	\begin{center}
		{\large\bf  Some Congruences Involving Binomial Coefficient and Fermat Quotient }
	\end{center}
	\vskip 2mm \centerline{Wei-Wei Qi}
	
	\begin{center}
		{\footnotesize MOE-LCSM, School of Mathematics and Statistics, Hunan Normal University, Hunan 410081, P.R. China\\[5pt]
			{\tt wwqi2022@foxmail.com} \\[10pt]
		}
	\end{center}
	
	\vskip 0.7cm \noindent{\bf Abstract.} In this paper,  we investigate some congruences involving sums of $\frac{d^{-k}{x\choose k}{x+k\choose k}}{{2k \choose k}}$, where $x$ be a $p$-adic integer, $k$ be a non-negative integer, and $d$ $(d\neq 0)$ be a rational number.

	\vskip 2mm
	\noindent{\it MR Subject Classifications}: 11A07, 11B65, 11B68, 05A10	
	
	\section{Introduction}

Recall that the harmonic numbers are defined by
\begin{align*}
H_n:=\sum_{k=0}^{n}\frac{1}{k}, \quad and \quad H_0:=0.
\end{align*}	
The Bernoulli numbers $\{B_n\}$ and Bernoulli polynomials $\{B_n(x)\}$ are defined by
\begin{align*}
 \sum_{k=0}^{n-1}{n\choose k}B_k=0 (n\geq 2)\quad and \quad B_n(x)=\sum_{k=0}^{n}{n\choose k}B_kx^{n-k}(n\geq 0).
\end{align*}
And the Euler numbers $\{E_n\}$ and Euler polynomials $\{E_n(x)\}$ are given by
\begin{align*}
\frac{2e^{t}}{1+e^{2t}}=\sum_{k=0}^{\infty}{n\choose k}E_k\frac{t^k}{k!} (\left|t\right|\textless \frac{\pi}{2})\quad and \quad \frac{2e^{xt}}{1+e^{t}}=\sum_{k=0}^{\infty}{n\choose k}E_k(x)\frac{t^k}{k!} (\left|t\right|\textless \pi).
\end{align*}
The reader is referred to \cite{w1-03} for more basic properties of the Bernoulli and Euler polynomials.

Congruences involving binomial coefficients are an interesting project, which have been studied widely by many authors. For more studies on binomial coefficients, see ( \cite{w1-03-2}--\cite{w1-08}, \cite{w1-09}--\cite{w1-012} and so on).	
In $2010$, L.-L. Zhao, H. Pan, and Z.-W. Sun \cite{w1-03-2} studied congruences for sums involving ${3k\choose k}$ and  proved that for any prime $p\geq 5$,	
\begin{align*}
	\sum_{k=0}^{p-1}{3k\choose k}2^k\equiv \frac{6(-1)^{(p-1)/2}-1}{5} \pmod{p}.
\end{align*}
Let $p$ be an odd prime and let $\mathbb{Z}_p$ stand for the ring of all $p$-adic integers. Z.-W. Sun \cite{w1-03-1} studied $\sum_{k=0}^{p-1}{3k \choose k+d}x^k\pmod{p}$. In particular, he \cite{w1-03-1} showed that
\begin{align*}
	\begin{aligned}
	\sum_{k=0}^{p-1}{3k \choose k+d}(\frac{4}{27})^k \equiv 
		\begin{cases}
			\frac{1}{9}  \pmod{p}  \quad&{if \quad  d=0},\\
			-\frac{16}{9}   \pmod{p} \quad&{if \quad  d=1},\\
			-\frac{4}{9}   \pmod{p} \quad&{if \quad  d=-1}.
		\end{cases}
	\end{aligned}
\end{align*}
Mattarei and Tauraso \cite{w1-03-3} deduced that for any prime $p\geq 3$
\begin{align*}
	\sum_{k=0}^{p-1}{2k\choose k}\equiv (\frac{p}{3})-\frac{1}{3}p^2B_{p-2}(\frac{1}{3}) \pmod{p^3},
\end{align*}
where $(\frac{\cdot}{\cdot})$ is the Jacobi symbol.	In $2015$, Kh. Hessami Pilehrood and T. Hessami Pilehrood \cite{w1-04} further studied congruences involving ${3k\choose k}$, ${4k\choose k}$ and the sequence (cf.\cite[A176898]{w1-05})
\begin{align*}
	S_k=\frac{{6k\choose 3k}{3k\choose k}}{2(2k+1){2k\choose k}},\quad k=0,1,2,\dots
\end{align*}
It is easy to see that
\begin{align*}
&{2k\choose k}=\frac{{-1/2\choose k}{-1/2+k\choose k}(-16)^k}{{2k\choose k}}, \quad\quad {3k\choose k}=\frac{{-1/3\choose k}{-1/3+k\choose k}(-27)^k}{{2k\choose k}}, \\
&{4k\choose 2k}=\frac{{-1/4\choose k}{-1/4+k\choose k}(-64)^k}{{2k\choose k}},\quad \frac{{6k\choose k}{3k\choose k}}{{2k\choose 2k}}=\frac{{-1/6\choose k}{-1/6+k\choose k}(-432)^k}{{2k\choose k}},
\end{align*}	
where 
\begin{align*}
	{x\choose k}=\frac{(x)(x-1)\cdots(x-k+1)}{k!} \quad (x\in \mathbb{B}, k\in \mathbb{N}={0,1,2,\cdots})
\end{align*}
are (generalized) binomial coefficients.

Recently, Wang and Han \cite[Theorem $1.2$]{w1-01} proved that for an odd prime $p$ and a $p$-adic number integer $x$,	
\begin{align*}
	\begin{aligned}
		\sum_{k=0}^{p-1}\frac{{x\choose k}{x+k\choose k}(-2)^k}{{2k\choose k}}&\equiv (-1)^{[(\langle x\rangle_p+1)/2]}\left(1+t-(-1)^{\langle x\rangle_p}(\frac{-1}{p})t\right)\\
		&-\frac{pt(t+1)}{2}\left(E_{p-2}(\frac{x+1}{2})+E_{p-2}(-\frac{x}{2})\right) \pmod{p^2},
	\end{aligned}
\end{align*}	
where  $x=\langle x\rangle_p+pt$, $t\in \mathbb{Z}_p$. Mao \cite[Theorem 1.2]{w1-02} also showed that for an odd prime $p$ and $p$-adic number integer $x$, then modulo $p^2$	
\begin{align*}
	\begin{aligned}
		\sum_{k=0}^{(p-1)/2}&\frac{{x\choose k}{x+k\choose k}(-2)^k}{{2k\choose k}}
		\equiv 
		\begin{cases}
			(-1)^{\lfloor\frac{\langle x \rangle_p+1}{2} \rfloor}+\frac{pt}{2}(\frac{-2}{p})E_{p-2}(\frac{2x+3}{4}) \quad&{if \quad \langle x \rangle_p\leq (p-1)/2},\\
			(\frac{-1}{p})(-1)^{\lfloor\frac{\langle x \rangle_p}{2} \rfloor}+\frac{p(1+t)}{2}(\frac{-2}{p})E_{p-2}(\frac{2x+3}{4}) \quad&{if \quad \langle x \rangle_p\textgreater (p-1)/2}.
		\end{cases} 
	\end{aligned}
\end{align*}	

Motivated by the above, in this paper, we further investigate supercongruence for sums involving ${x\choose k}{x+k\choose k}/{2k\choose k}$.	Throughout, for any prime $p$ and $x\in \mathbb{Z}_p$, we always use $\langle x\rangle_p$  denote the least nonnegative residue of $x$ modulo $p$. Write $x=\langle x\rangle_p+mp$, $m\in \mathbb{Z}_p$. The Fermat quotient of an integer $a$ with respect to an odd prime $p$ is given by 
\begin{align*}
	q_p(a)=\frac{a^{p-1}-1}{p}.
\end{align*}	

	\begin{thm}
	Let $p\geq 3$ be a prime, $d$ ($d\neq 0$) be a rational number, $x$ be a $p$-adic integer and  $m:=(x-\langle x\rangle_p)/p$, $\langle x\rangle_p \in \{0,1,2,\dots p-1\}$. Then modulo $p^4$ 
		\begin{align}
			\begin{aligned}
			\sum_{k=0}^{p-1}&\left(x^2+x-(1+4d)k^2-(1-2d)k\right)\frac{d^{-k}{x\choose k}{x+k\choose k}}{{2k\choose k}}\\
			&\equiv d^{1-p}m(1+m)\left(2p^2-p-2p^2(1-2p)H_{\langle x\rangle_p}+2p^3mH_{\langle x\rangle_p}^{(2)}-2p^3H_{\langle x\rangle_p}^2\right). \label{th-1}
				\end{aligned}
		\end{align}
	\end{thm}			
	
Taking $x=-1/2$, $-1/3$, $-1/4$ and $-1/6$ in \eqref{th-1}, we have the following consequences.	
\begin{cor}	Let $p\textgreater 3$ be a prime and $d$ ($d\neq 0$) be a rational number. Then modulo $p^4$
\begin{align}
	\begin{aligned}
	\sum_{k=0}^{p-1}\left(1+(4-d)k^2+\frac{8+d}{2}k\right)\frac{{2k\choose k}}{d^k}\equiv (\frac{16}{d})^{p-1}\left(p(2p-1)(1-4pq_p(2))-10p^3q_p(2)^2\right),\label{c0-1}
	\end{aligned}
\end{align}
\begin{align}
	\begin{aligned}
		\sum_{k=0}^{p-1}\left(1+\frac{27-4d}{6}k^2+\frac{27+2d}{6}k\right)\frac{{3k\choose k}}{d^k}\equiv (-\frac{27}{d})^{p-1}p\left((2p-1)(1-3pq_p(3))-6p^2q_p(3)^2\right), \label{c0-2}
	\end{aligned}
\end{align}
\begin{align}
	\begin{aligned}
		\sum_{k=0}^{p-1}\left(1+\frac{16-d}{3}k^2+\frac{32+d}{6}k\right)\frac{{4k\choose 2k}}{d^k}\equiv (-\frac{64}{d})^{p-1}p\left((2p-1)(1-6pq_p(2))-21p^2q_p(2)^2\right), \label{c0-3}
	\end{aligned}
\end{align}
\begin{align}
	\begin{aligned}
		\sum_{k=0}^{p-1}&\left(1+\frac{108-d}{15}k^2+\frac{216+d}{30}k\right)\frac{{4k\choose 2k}{6k\choose 3k}}{d^k{2k\choose k}}\\
		&\equiv (-\frac{432}{d})^{p-1}p\left((2p-1)(1-4pq_p(2)-3pq_p(3))-2p^2(5q_p(2)^2+6q_p(2)q_p(3)+3q_p(3)^2)\right). \label{c0-4}
	\end{aligned}
\end{align}	
\end{cor}

	\begin{thm}
Let $p\geq 3$ be a prime, $d$ ($d\neq 0$) be a rational number, $x$ be a $p$-adic integer and  $m:=(x-\langle x\rangle_p)/p$, $\langle x\rangle_p \in \{0,1,2,\dots p-1\}$.
	 If $\langle x\rangle_p \textless (p-1)/2$, then modulo $p^4$ 
	\begin{align}
		\begin{aligned}
			\sum_{k=0}^{\frac{p-1}{2}}&\left(x^2+x-(1+4d)k^2-(1-2d)k\right)\frac{d^{-k}{x\choose k}{x+k\choose k}}{{2k\choose k}}\\
			&\equiv \frac{(-1)^{\langle x\rangle_p}m(-\frac{1}{d})^{(p-1)/2}}{2{p-1\choose \langle x\rangle_p+(p-1)/2}}
			\left(p+2\langle x\rangle_p p+(1+2m)p^2+(1+2\langle x\rangle_p)mp^3(H_{\langle x\rangle_p}^{(2)}-4H_{2\langle x\rangle_p}^{(2)})\right).  \label{th-2-1}
		\end{aligned}
	\end{align}
If $\langle x\rangle_p = (p-1)/2$, then modulo $p^5$
\begin{align}
	\begin{aligned}
		\sum_{k=0}^{\frac{p-1}{2}}&\left(x^2+x-(1+4d)k^2-(1-2d)k\right)\frac{d^{-k}{x\choose k}{x+k\choose k}}{{2k\choose k}}
		\equiv  \frac{m(1+m)p^2}{4^{p-1}d^{\frac{p-1}{2}}}\left(1+2pq_p(2)+p^2q_p(2)^2\right) . \label{th-2-2}
	\end{aligned}
\end{align}
If $\langle x\rangle_p \textgreater (p-1)/2$, then modulo $p^2$
\begin{align}
	\begin{aligned}
		\sum_{k=0}^{\frac{p-1}{2}}&\left(x^2+x-(1+4d)k^2-(1-2d)k\right)\frac{d^{-k}{x\choose k}{x+k\choose k}}{{2k\choose k}}
		\equiv\frac{(1+m)p}{d^{\frac{p-1}{2}}} {\langle x\rangle_p+\frac{p+1}{2}\choose p+1}. \label{th-2-3}
	\end{aligned}
\end{align}
\end{thm}			
	
Putting $x=-1/2$, $-1/3$, $-1/4$ and $-1/6$ in Theorem $1.3$, we have the following results.
\begin{cor}	Let $p\textgreater 3$ be a prime and $d$ ($d\neq 0$) be a rational number. Then 
	\begin{align}
		\begin{aligned}
		\sum_{k=0}^{(p-1)/2}\left(1+(4-d)k^2+\frac{8+d}{2}k\right)\frac{{2k\choose k}}{d^k}\equiv(-\frac{1}{d})^{(p-1)/2}p^2 \left(1+2pq_p(2)+p^2q_p(2)^2\right) \pmod{p^5},\label{c0-5}
		\end{aligned}
	\end{align}
\begin{align}
		\begin{aligned}
\sum_{k=0}^{(p-1)/2}&\left(1+\frac{27-4d}{6}k^2+\frac{27+2d}{6}k\right)\frac{{3k\choose k}}{d^k}\\
& \equiv 
	\begin{cases}
		\frac{(-1)^{(p-1)/3}(\frac{27}{d})^{(p-1)/2}}{4{p-1\choose (5p-5)/6}}\left(p+3p^2-\frac{p^3}{6}(\frac{p}{3})B_{p-2}(\frac{1}{6})\right) \pmod{p^4}\quad&{if \quad p\equiv 1\pmod 3},\\
		-\frac{3}{2}(-\frac{27}{d})^{(p-1)/2}p{(7p+1)/6 \choose p+1}	 \pmod{p^2} \quad&{if\quad  p\equiv 2\pmod 3},\label{c0-6}
	\end{cases} 
	\end{aligned}
\end{align}	
\begin{align}
	\begin{aligned}
		\sum_{k=0}^{(p-1)/2}&\left(1+\frac{16-d}{3}k^2+\frac{32+d}{6}k\right)\frac{{4k\choose 2k}}{d^k}\\
	& \equiv 
	\begin{cases}
		\frac{(-1)^{(p-1)/4}(\frac{64}{d})^{(p-1)/2}}{3{p-1\choose (3p-3)/4}}\left(p+2p^2-(-1)^{(p-1)/2}p^3E_{p-3}\right) \pmod{p^4}\quad&{if \quad p\equiv 1\pmod 4},\\
		-\frac{4}{3}(-\frac{64}{d})^{(p-1)/2}p{(5p+1)/4 \choose p+1}	 \pmod{p^2} \quad&{if\quad  p\equiv 3\pmod 4},\label{c0-7}
	\end{cases} 
	\end{aligned}
\end{align}	
\begin{align}
		\begin{aligned}
	\sum_{k=0}^{(p-1)/2}&\left(1+\frac{108-d}{15}k^2+\frac{216+d}{30}k\right)\frac{{4k\choose 2k}{6k\choose 3k}}{d^k{2k\choose k}}\\
	& \equiv 
	\begin{cases}
		\frac{(-1)^{(p-1)/6}(\frac{432}{d})^{(p-1)/2}}{5{p-1\choose (2p-2)/3}}\left(2p+3p^2-\frac{p^3}{30}(\frac{p}{3})B_{p-2}(\frac{1}{6})\right) \pmod{p^4}\quad&{if \quad p\equiv 1\pmod 6},\\
		-\frac{6}{5}(-\frac{432}{d})^{(p-1)/2}p{(4p+1)/3 \choose p+1}	 \pmod{p^2} \quad&{if\quad  p\equiv 5\pmod 6}.\label{c0-8}
	\end{cases}
 	\end{aligned}
\end{align}	
\end{cor}

	\begin{thm}
	(i).Let $n$ be a positive integer, $d$ ($d\neq 0$) be a rational number, $x$ be a $p$-adic integer. Then
	\begin{align}
		\begin{aligned}
			\sum_{k=0}^{n-1}&\left(x^2+x-(1+4d)k^2-(1-2d)k\right)\frac{d^{-k}{x\choose k}{x+k\choose k}}{{2k\choose k}}\\
			&=-2nd+n\sum_{k=0}^{n-1}\left(x^2+x+2d-(1+4d)k^2-(1+2d)k\right)\frac{d^{-k}{x\choose k}{x+k\choose k}}{(k+1){2k\choose k}}\\
			&=(2n-1)\sum_{k=0}^{n-1}\left(x^2+x-(1+4d)k^2-(1+2d)k\right)\frac{d^{-k}{x\choose k}{x+k\choose k}}{(2k+1){2k\choose k}}. \label{th-1-0}
		\end{aligned}
	\end{align}	
	(ii).Let $n$ be a positive odd number, $d$ ($d\neq 0$) be a rational number, $x$ be a $p$-adic integer. Then
\begin{align}
	\begin{aligned}
		\sum_{k=0}^{\frac{n-1}{2}}&\left(x^2+x-(1+4d)k^2-(1-2d)k\right)\frac{d^{-k}{x\choose k}{x+k\choose k}}{{2k\choose k}}\\
		&=-d(n+1)+\frac{n+1}{2}\sum_{k=0}^{\frac{n-1}{2}}\left(x^2+x+2d-(1+4d)k^2-(1+2d)k\right)\frac{d^{-k}{x\choose k}{x+k\choose k}}{(k+1){2k\choose k}}\\
		&=n\sum_{k=0}^{\frac{n-1}{2}}\left(x^2+x-(1+4d)k^2-(1+2d)k\right)\frac{d^{-k}{x\choose k}{x+k\choose k}}{(2k+1){2k\choose k}}. \label{th-2-0}
	\end{aligned}
\end{align}	
\end{thm}		

\textit{Remark}. For $p$-adic integer $x$, rational number $d$ ($d\neq 0$). Combining Theorem $1.1$, Theorem $1.3$ and Theorem $1.5$, we can obtain similar congruences  for  the sum of polynomial families involving  $\frac{d^{-k}{x\choose k}{x+k\choose k}}{(2k+1){2k\choose k}}$ and $\frac{d^{-k}{x\choose k}{x+k\choose k}}{(k+1){2k\choose k}}$.

The rest of the paper is organized as follows. We shall prove Theorems  and Corollaries in the section $2$ and section $3$, respectively.

\section{Proofs of the Theorems}
For each positive integer $n$ and $r$, letting	
\begin{align*}
	H_n^{(r)}=\sum_{k=1}^{n}\frac{1}{k^r}, \quad and \quad H(1,1;n)=\sum_{k=1}^{n}\frac{H_{k-1}}{k}.
\end{align*}

In order to achieve the proofs of Theorems, we need the following Lemma:	
	\begin{lem}
		Let $d$ $(d\neq 0)$ be a rational number, $x$ be a $p$-adic integer,  $n$ be non-negative integer. Then
		\begin{align}
			\sum_{k=0}^{n-1}\left(x^2+x-(1+4d)k^2-(1-2d)k\right)\frac{d^{-k}{x\choose k}{x+k\choose k}}{{2k\choose k}}=\frac{2n(2n-1){x\choose n}{x+n\choose n}}{d^{n-1}{2n\choose n}}, \label{lm-1}
		\end{align}
\begin{align}
	\sum_{k=0}^{n-1}\left(x^2+x-(1+4d)k^2-(1+2d)k\right)\frac{d^{-k}{x\choose k}{x+k\choose k}}{(2k+1){2k\choose k}}=\frac{2n{x\choose n}{x+n\choose n}}{d^{n-1}{2n\choose n}}, \label{lm-1-0}
\end{align}
and
\begin{align}
	\sum_{k=0}^{n-1}\left(x^2+x+2d-(1+4d)k^2-(1+2d)k\right)\frac{d^{-k}{x\choose k}{x+k\choose k}}{(k+1){2k\choose k}}=2d+\frac{2(2n-1){x\choose n}{x+n\choose n}}{d^{n-1}{2n\choose n}}. \label{lm-1-0-0}
\end{align}
	\end{lem}
	\pf. Noting that	
	\begin{align*}	
	 (2k+1)(2k+2)\frac{{x\choose k+1}{x+k+1\choose k+1}}{{2k+2\choose k+1}}=(x-k)(x+k+1)\frac{{x\choose k}{x+k\choose k}}{{2k\choose k}}.	
	\end{align*}
For real number $d$, we have
\begin{align}	
(2k+1)(2k+2)d^{n-1-k}\frac{{x\choose k+1}{x+k+1\choose k+1}} {{2k+2\choose k+1}}=(x-k)(x+k+1)d^{n-1-k}\frac{{x\choose k}{x+k\choose k}}{{2k\choose k}}. \label{lm-1-01}
\end{align}
Summing both sides of the above over $k$ from $0$ to $n-1$ gives
	\begin{align}
	\begin{aligned}	
		\sum_{k=0}^{n-1}(2k+1)(2k+2)d^{n-1-k}\frac{{x\choose k+1}{x+k+1\choose k+1}} {{2k+2\choose k+1}}=\sum_{k=0}^{n-1}(x-k)(x+k+1)d^{n-1-k}\frac{{x\choose k}{x+k\choose k}}{{2k\choose k}},\label{lm-1-02}
	\end{aligned}	
\end{align}
the left hand-side could be changed to
\begin{align}
	\begin{aligned}	
	\sum_{k=0}^{n-1}&(2k+1)(2k+2)d^{n-1-k}\frac{{x\choose k+1}{x+k+1\choose k+1}} {{2k+2\choose k+1}}\\
	&=\sum_{k=0}^{n-1}(2k-1)(2k)d^{n-k}\frac{{x\choose k}{x+k\choose k}} {{2k\choose k}}+(2n-1)(2n)\frac{{x\choose n}{x+n\choose n}} {{2n\choose n}}.\label{lm-1-03}
	\end{aligned}	
\end{align}
Combining \eqref{lm-1-02} and \eqref{lm-1-03} gives \eqref{lm-1}.\\
By \eqref{lm-1-01}, we have
\begin{align*}
	\begin{aligned}	
		\sum_{k=0}^{n-1}&(x-k)(x+k+1)d^{n-1-k}\frac{{x\choose k}{x+k\choose k}}{(2k+1){2k\choose k}}\\
		&=\sum_{k=0}^{n-1}(2k+2)d^{n-1-k}\frac{{x\choose k+1}{x+k+1\choose k+1}} {{2k+2\choose k+1}}\\
		&=\sum_{k=0}^{n-1}2kd^{n-k}\frac{{x\choose k}{x+k\choose k}} {{2k\choose k}}+2n\frac{{x\choose n}{x+n\choose n}} {{2n\choose n}}.
	\end{aligned}	
\end{align*}
From the above, we get \eqref{lm-1-0}. The proof of \eqref{lm-1-0-0} is similar to that of \eqref{lm-1-0}, so we omit the details. This confirms Lemma $2.1$. \qed \\

\textit{Proof of \eqref{th-1}}.
 Considering that
\begin{align*}
	H(1,1;\langle x\rangle_p)=\sum_{k=1}^{\langle x\rangle_p}\frac{H_{k-1}}{k}=\frac{1}{2}\left(H_{\langle x\rangle_p}^2-H_{\langle x\rangle_p}^{(2)}\right),
\end{align*}
and for $p\geq 3$ we have
\begin{align*}
	H_{p-1}\equiv 0 \pmod{p^2}, \quad \quad \sum_{k=1}^{\langle x\rangle_p}\frac{1}{p-k}\equiv -\sum_{k=1}^{\langle x\rangle_p}\left(\frac{1}{k}+\frac{p}{k^2}\right) \pmod{p^2}.
\end{align*}
It following that
\begin{align}
	\begin{aligned}	
		&{pm+\langle x\rangle_p \choose p}
		{p(m+1)+\langle x\rangle_p\choose p}\\
		&=\frac{(pm+p+\langle x\rangle_p)(pm+p+\langle x\rangle_p-1)\dots(pm-(p-\langle x\rangle_p-1))}{p!p!}\\
		&\equiv \frac{m(m+1)(-1)^{\langle x\rangle_p}}{{p-1\choose \langle x\rangle_p}}\left(1+(pm+p)H_{\langle x\rangle_p}+(pm+p)^2H(1,1;\langle x\rangle_p)\right)\\
		&\times \left(1-pmH_{p-1-\langle x\rangle_p}+p^2m^2H(1,1;p-1-\langle x\rangle_p)\right)\\
		&\equiv m(m+1)\left(1+2pH_{\langle x\rangle_p}-2mp^2H_{\langle x\rangle_p}^{(2)}+2p^2H_{\langle x\rangle_p}^2\right) \pmod{p^3}. \label{l-1}
	\end{aligned}	
\end{align}

 In $1862$, J. Wolstenholme \cite{w1} proved that for $p\geq 3$,	
\begin{align}
	\begin{aligned}
	{2p-1\choose p-1}\equiv 1 \pmod{p^3}. \label{l-2}
	\end{aligned}
\end{align}			
Setting $n=p$ in \eqref{lm-1}, with the help of \eqref{l-1} and \eqref{l-2} gives	
\begin{align*}
	\begin{aligned}
			&\sum_{k=0}^{p-1}\left(x^2+x-(1+4d)k^2-(1-2d)k\right)\frac{d^{-k}{x\choose k}{x+k\choose k}}{{2k\choose k}}\\
			&=\frac{2p(2p-1)}{d^{p-1}{2p\choose p}}{x\choose p}{x+p\choose p}
		\\
		&\equiv d^{1-p}(2p-1)p(m^2+m)\left(1+2pH_{\langle x\rangle_p}-2p^2mH_{\langle x\rangle_p}^{(2)}+2p^2H_{\langle x\rangle_p}^2\right) \pmod{p^4},
	\end{aligned}
\end{align*}	
as desired.	This, the proof of Theorem $1.1$ is complete.\\	

 Furthermore, setting $n=\frac{p+1}{2}$ in \eqref{lm-1} and simplifying, we get that for prime $p\geq 3$,
\begin{align}
	\begin{aligned}
		\sum_{k=0}^{\frac{p-1}{2}}&\left(x^2+x-(1+4d)k^2-(1-2d)k\right)\frac{d^{-k}{x\choose k}{x+k\choose k}}{{2k\choose k}}
		=\frac{6pd^{\frac{1-p}{2}}}{p^2-1}{x\choose \frac{p+1}{2}}{x+\frac{p+1}{2}\choose \frac{p+1}{2}}/{p-1\choose \frac{p-1}{2}}. \label{l-3}
	\end{aligned}
\end{align}	
Next, we divide into three cases based on value of $\langle x\rangle_p$.

\textbf{Case 1.} If $\langle x\rangle_p \textless \frac{p-1}{2}$.\\

In \cite{w2}, L. Carlitz showed that for prime $p\geq3$
\begin{align}
	\begin{aligned}
		{p-1\choose \frac{p-1}{2}}\equiv (-1)^{\frac{p-1}{2}}\left(4^{p-1}+\frac{p^3}{12}B_{p-1}\right) \pmod{p^4}.\label{l-4}
	\end{aligned}
\end{align}
Hence, we have
\begin{align*}
	\begin{aligned}
		\frac{{\langle x\rangle_p+(p+1)/2\choose \langle x\rangle_p}}{{(p-1)/2 \choose \langle x\rangle_p}}&=\frac{{p-1\choose (p-1)/2}}{{p-1\choose \langle x\rangle_p+(p-1)/2}}\frac{1+p+2\langle x\rangle_p}{1+p}
		\equiv \frac{(-1)^{\frac{p-1}{2}}4^{p-1}}{{p-1\choose \langle x\rangle_p+(p-1)/2}}\left(1+2\langle x\rangle_p(1-p+p^2)\right)\pmod{p^3}. 
	\end{aligned}
\end{align*}		
In addition, through calculation we obtain  	
\begin{align*}
	\begin{aligned}
		H_{\langle x\rangle_p+(p+1)/2}-H_{(p-1)/2-\langle x\rangle_p}&=\sum_{k=1}^{\langle x\rangle_p}\frac{2}{p+1-2k}+\sum_{k=1}^{\langle x\rangle_p}\frac{2}{p+1+2k}+\frac{2}{1+p}\\
		&\equiv p\sum_{k=1}^{\langle x\rangle_p}\frac{1}{k^2}-4p\sum_{k=1}^{2\langle x\rangle_p}\frac{1}{k^2}+\frac{4\langle x\rangle_p-2p+2}{(1+2\langle x\rangle_p)^2} \pmod{p^2}.
	\end{aligned}
\end{align*}		
It is not difficult for us to compute that
\begin{align*}
	\begin{aligned}
		H(1,1;\langle x\rangle_p+(p+1)/2)
		-H_{\langle x\rangle_p+(p+1)/2}H_{(p-1)/2-\langle x\rangle_p}+H(1,1;(p-1)/2-\langle x\rangle_p)\equiv 0 \pmod{p}.
	\end{aligned}
\end{align*}		
It following that
\begin{align}
	\begin{aligned}
		&{pm+\langle x\rangle_p \choose \frac{p+1}{2}}{pm+\langle x\rangle_p +\frac{p+1}{2}\choose \frac{p+1}{2}}\\
		&=\frac{(pm+\langle x\rangle_p+(p+1)/2)\dots pm(pm-1)\dots((pm+\langle x\rangle_p-(p-1)/2))}{(p+1)/2!(p+1)/2!}\\
		&\equiv \frac{2{\langle x\rangle_p+(p+1)/2 \choose (p+1)/2}}{{(p-1)/2\choose \langle x\rangle_p}(1+p)}(-1)^{(p-1)/2-\langle x\rangle_p}pm\\
		&\times((1+mp(H_{\langle x\rangle_p+(p+1)/2}-H_{(p-1)/2-\langle x\rangle_p})+m^2p^2(H(1,1;\langle x\rangle_p+(p+1)/2)\\
		&-H_{\langle x\rangle_p+(p+1)/2}H_{(p-1)/2-\langle x\rangle_p}+H(1,1;(p-1)/2-\langle x\rangle_p)))\\
		&\equiv \frac{(-1)^{\langle x\rangle_p } 4^{p-1}2pm}{{p-1\choose \langle x\rangle_p+(p-1)/2}(1+p)}\\
		&\times\left(1+2\langle x\rangle_p+p(2m-2\langle x\rangle_p-2pm+2p\langle x\rangle_p)+mp^2(1+2\langle x\rangle_p)(H_{\langle x\rangle_p}^{(2)}-4H_{2\langle x\rangle_p}^{(2)})\right) \pmod{p^4}. \label{l-5}
	\end{aligned}
\end{align}		
Then,  substituting \eqref{l-4} and \eqref{l-5}	into \eqref{l-3} and simplifying, we arrive at \eqref{th-2-1}.\\

\textbf{Case 2.} If $\langle x\rangle_p = \frac{p-1}{2}$.

In view of \cite[Lemma 2.2]{w3}: For prime $p\geq 3$	
\begin{align}
	\sum_{k=1}^{(p-1)/2}\frac{1}{k}\equiv -2q_p(2)+pq_p(2)^2 \pmod{p^2}. \label{l-5-1}
\end{align}		
Meanwhile, we have	
\begin{align}
	\sum_{k=1}^{(p-1)/2}\frac{1}{k^2}\equiv \frac{1}{2}\sum_{k=1}^{p-1}\frac{1}{k^2}\equiv 0 \pmod{p}. \label{l-5-2}
\end{align}	
Therefore, by \eqref{l-5-1} and \eqref{l-5-2} yields	
\begin{align}
	\begin{aligned}
		&{pm+(p-1)/2 \choose (p+1)/2 }{pm+p\choose (p+1)/2}\\
		&=p^2m(1+m)\frac{(pm+(p-1)/2)\dots (pm+1)(pm+p-1)\dots((pm+p-(p-1)/2))}{(p+1)/2!(p+1)/2!}\\
		&\equiv \frac{(-1)^{\frac{p-1}{2}}4m(1+m)p^2}{(1+p)^2}\left(1-pH_{(p-1)/2}+\frac{p^2}{2}H_{(p-1)/2}^2-\frac{(1+2m+2m^2)p^2}{2}H_{(p-1)/2}^{(2)}\right)\\
		&\equiv \frac{4(-1)^{\frac{p-1}{2}}m(1+m)p^2}{(1+p)^2}\left(1+2pq_p(2)+p^2q_p(2)^2\right) \pmod{p^5}. \label{l-6}
	\end{aligned}
\end{align}			
Combining \eqref{l-3}, \eqref{l-4} and \eqref{l-6}, we obtain	
\begin{align*}
	\begin{aligned}
		\sum_{k=0}^{\frac{p-1}{2}}&\left(x^2+x-(1+4d)k^2-(1-2d)k\right)\frac{d^{-k}{x\choose k}{x+k\choose k}}{{2k\choose k}}\\
		&\equiv  \frac{m(1+m)p^2}{4^{p-1}d^{\frac{p-1}{2}}}\left(1+2pq_p(2)+p^2q_p(2)^2\right) \pmod{p^5},
	\end{aligned}
\end{align*}	
as claimed.	\\

\textbf{Case 3.} If $\langle x\rangle_p \textgreater \frac{p-1}{2}$.\\	

Considering that
\begin{align*}
	\begin{aligned}
		\sum_{k=1}^{\langle x\rangle_p+(p+1)/2}\frac{1}{k}-\sum_{k=1}^{\langle x\rangle_p-(p-1)/2}\frac{1}{k}&=H_{(p-1)/2}+\sum_{k=1}^{\langle x\rangle_p+1}\frac{2}{2k+p-1}-\sum_{k=1+(p-1)/2}^{\langle x\rangle_p}\frac{2}{2k-p+1}\\
		&\equiv \sum_{k=1}^{(p-1)/2}\frac{1}{k}+\sum_{k=1}^{(p-1)/2}\frac{2}{2k-1}+\frac{1}{p}\\
		&\equiv \frac{1}{p} \pmod{p},
	\end{aligned}
\end{align*}	
where we also used the fact $H_{p-1}\equiv 0 \pmod{p}$ in the last step.\\
Moreover, utilize \eqref{l-4} we have
\begin{align*}
	\begin{aligned}
		{\langle x\rangle_p \choose \frac{p+1}{2}}{\langle x\rangle_p+\frac{p+1}{2}\choose \frac{p+1}{2}}=\frac{4p}{p+1}{p-1\choose \frac{p-1}{2}}{\langle x\rangle_p+\frac{p+1}{2}\choose p+1}\equiv (-1)^{\frac{p-1}{2}}4p{\langle x\rangle_p+\frac{p+1}{2}\choose p+1} \pmod{p^2}.
	\end{aligned}
\end{align*}	
Hence, by the above congruences gives
\begin{align}
	\begin{aligned}
		&{mp+\langle x\rangle_p \choose \frac{p+1}{2}}{mp+\frac{p+1}{2}+\langle x\rangle_p \choose \frac{p+1}{2}}\\
		&=
		\frac{(mp+(p+1)/2+\langle x\rangle_p)(mp+(p-1)/2+\langle x\rangle_p)\dots(mp+\langle x\rangle_p-(p-1)/2)}{((p-1)/2)!((p-1)/2)!}\\
		&\equiv {\langle x\rangle_p \choose \frac{p+1}{2}}{\langle x\rangle_p+\frac{p+1}{2}\choose \frac{p+1}{2}}\left(1+mp(\sum_{k=1}^{\langle x\rangle_p+(p+1)/2}\frac{1}{k}-\sum_{k=1}^{\langle x\rangle_p-(p-1)/2}\frac{1}{k})\right)\\
		&\equiv (-1)^{\frac{p-1}{2}}4p(1+m){\langle x\rangle_p+\frac{p+1}{2}\choose p+1} \pmod{p^2}.\label{l-7}
	\end{aligned}
\end{align}		
Finally, substituting   \eqref{l-4} and \eqref{l-7} into \eqref{l-3} and simplifying, we obtain	desired result. In view of the above, the proof of Theorem $1.3$ is now complete.

\textit{Proof of Theorem 1.5}. In view of \eqref{lm-1} and \eqref{lm-1-0}, we immediately get \eqref{th-1-0} and \eqref{th-2-0}.
Thus, the proofs are completed. \qed
	
\section{Proofs of the corollaries}
In order to show the proofs of Corollaries, we need the following Lemma:
\begin{lem}\cite[Corollaries 3.3, 3.7 and Theorem 3.9]{w1-0} For prime $p\textgreater 3$,
	\begin{align}
		\begin{aligned}
			&H_{\lfloor \frac{p}{3} \rfloor}\equiv -\frac{3}{2}q_p(3)+\frac{3}{4}pq_p(3)^2-\frac{p}{30}(\frac{p}{3})B_{p-2}(\frac{1}{6})\pmod{p^2},\\
			&H_{\lfloor \frac{p}{4} \rfloor}\equiv -3q_p(2)+\frac{3}{2}pq_p(2)^2-(-1)^{\frac{p-1}{2}}pE_{p-3}\pmod{p^2},\\
			&H_{\lfloor \frac{p}{6} \rfloor}\equiv -2q_p(2)-\frac{3}{2}q_p(3)+pq_p(2)^2+\frac{3}{4}pq_p(3)^2-\frac{p}{12}(\frac{p}{3})B_{p-2}(\frac{1}{6})\pmod{p^2},\\
			&H_{\lfloor \frac{p}{3} \rfloor}^{(2)}\equiv \frac{1}{10}(\frac{p}{3})B_{p-2}(\frac{1}{6})\pmod{p},\\
			&H_{\lfloor \frac{p}{4} \rfloor}^{(2)}\equiv (-1)^{\frac{p-1}{2}}4E_{p-3}\pmod{p},\\
			&H_{\lfloor \frac{p}{6} \rfloor}^{(2)}\equiv \frac{1}{2}(\frac{p}{3})B_{p-2}(\frac{1}{6})\pmod{p}.  \label{lm-2}
		\end{aligned}
	\end{align}
\end{lem}

\textit{Proof of Corollary 1.3}.

When $x=-1/3$, $-1/4$ and $-1/6$, we have
\begin{align}
	\begin{cases}
		\langle -\frac{1}{3}\rangle_p=(p-1)/3, \quad  m=\frac{-1/3-\langle -\frac{1}{3}\rangle_p}{p}=-1/3 \quad&{if \quad p\equiv 1\pmod 3},\\
		\langle -\frac{1}{3}\rangle_p=(2p-1)/3, \quad  m=\frac{-1/3-\langle -\frac{1}{3}\rangle_p}{p}=-2/3	 \quad&{if\quad  p\equiv 2\pmod 3},	\label{l-2-1}
	\end{cases}
\end{align}
\begin{align*}
	\begin{cases}
		\langle -\frac{1}{4}\rangle_p=(p-1)/4, \quad  m=\frac{-1/3-\langle -\frac{1}{4}\rangle_p}{p}=-1/4 \quad&{if \quad p\equiv 1\pmod 4},\\
		\langle -\frac{1}{4}\rangle_p=(3p-1)/4, \quad  m=\frac{-1/4-\langle -\frac{1}{4}\rangle_p}{p}=-3/4	 \quad&{if\quad  p\equiv 3\pmod 4},	
	\end{cases}
\end{align*}
and
\begin{align*}
	\begin{cases}
		\langle -\frac{1}{6}\rangle_p=(p-1)/6, \quad  m=\frac{-1/6-\langle -\frac{1}{6}\rangle_p}{p}=-1/6 \quad&{if \quad p\equiv 1\pmod 6},\\
		\langle -\frac{1}{6}\rangle_p=(5p-1)/6, \quad  m=\frac{-1/3-\langle -\frac{1}{3}\rangle_p}{p}=-5/6	 \quad&{if\quad  p\equiv 5\pmod 6}.	
	\end{cases}
\end{align*}
Meanwhile, we have
\begin{align*}
	\langle -\frac{1}{2}\rangle_p=(p-1)/2, \quad  m=\frac{-1/2-\langle -\frac{1}{2}\rangle_p}{p}=-1/2,	
\end{align*}
Hence, letting $d\rightarrow-\frac{d}{16}$ in \eqref{th-1}, substituting the above, \eqref{l-5-1} and  \eqref{l-5-2} into \eqref{th-1} and with some necessary calculation, we get \eqref{c0-1}. 

Moreover, for $0\leq k \leq p-1$ we have
\begin{align*}
	H_{p-1-k}=\sum_{i=k+1}^{p-1}\frac{p+i}{p^2-i^2}\equiv -\sum_{i=k+1}^{p-1}\frac{p+i}{i^2} \equiv H_k+pH_k^{(2)}\pmod{p^2}. 
\end{align*}
Also, we have
\begin{align*}
	H_{p-1-k}^{(2)}=\sum_{i=k+1}^{p-1}\frac{1}{(p-i)^2}\equiv -H_k{(2)}\pmod{p}.
\end{align*}
Thus, by the above and \eqref{lm-2} we obtain
\begin{align}
	H_{(2p-1)/3}^{(2)}\equiv -H_{\lfloor (p-2)/3\rfloor }\equiv -\frac{1}{10}(\frac{p}{3})B_{p-2}(\frac{1}{6})\pmod{p}, \label{l-2-2}
\end{align}
and
\begin{align}
	\begin{aligned}
		H_{(2p-1)/3}&\equiv H_{\lfloor (p-2)/3\rfloor}+pH_{\lfloor (p-2)/3\rfloor}^{(2)}\\
		&\equiv-\frac{3}{2}q_p(3)+\frac{3}{4}pq_p(3)^2+\frac{p}{15}(\frac{p}{3})B_{p-2}(\frac{1}{6})\pmod{p^2}. \label{l-2-3}
	\end{aligned}
\end{align}
Letting $d\rightarrow-\frac{d}{27}$, substituting \eqref{l-2-1} into \eqref{th-1} and simplifying, then we using \eqref{lm-2}, \eqref{l-2-2},  \eqref{l-2-3} to obtain \eqref{c0-2}. The proof of \eqref{c0-3} and \eqref{c0-4} can be proceed as the \eqref{c0-2}, we omit them. This, the proof is complete. \qed

\textit{Proof of Corollary 1.4}.

If $x=-1/2$,  letting $d\rightarrow-\frac{d}{16}$ in \eqref{th-2-2} and simplifying, we immediately get claimed result.	
	
If 	$x=-1/3$, then for $p\equiv 1\pmod{3}$,	we have $\langle -\frac{1}{3}\rangle_p=(p-1)/3,   m=\frac{-1/3-\langle -\frac{1}{3}\rangle_p}{p}=-1/3$. Letting $d\rightarrow-\frac{d}{27}$ in \eqref{th-2-1}, and using \eqref{lm-2} yields
\begin{align}
	\sum_{k=0}^{(p-1)/2}&\left(1+\frac{27-4d}{6}k^2+\frac{27+2d}{6}k\right)\frac{{3k\choose k}}{d^k}\\
	& \equiv 
		\frac{(-1)^{(p-1)/3}(\frac{27}{d})^{(p-1)/2}}{4{p-1\choose (5p-5)/6}}\left(p+3p^2-\frac{p^3}{6}(\frac{p}{3})B_{p-2}(\frac{1}{6})\right) \pmod{p^4}.
\end{align}		
For $p\equiv 2\pmod{3}$, $\langle -\frac{1}{3}\rangle_p=(2p-1)/3,   m=-2/3$. Letting $d\rightarrow-\frac{d}{27}$ in \eqref{th-2-3}, we get
\begin{align}
	\sum_{k=0}^{(p-1)/2}&\left(1+\frac{27-4d}{6}k^2+\frac{27+2d}{6}k\right)\frac{{3k\choose k}}{d^k}
	\equiv	-\frac{3}{2}(-\frac{27}{d})^{(p-1)/2}p{(7p+1)/6 \choose p+1}	 \pmod{p^2}, 
\end{align}		
as desired. The proof of \eqref{c0-7} and \eqref{c0-8} are obtained similar to the proof of \eqref{c0-6}, so we overleap them. Thus, the proof is finished.\qed

\end{document}